# Controlled Islanding via Weak Submodularity

Zhipeng Liu, *Student Member, IEEE,* Andrew Clark, *Member, IEEE,* Linda Bushnell, *Fellow, IEEE,* Daniel Kirschen, *Fellow, IEEE,* and Radha Poovendran, *Fellow, IEEE*

*Abstract*—Cascading failures typically occur following a large disturbance in power systems, such as tripping of a generating unit or a transmission line. Such failures can propagate and destabilize the entire power system, potentially leading to widespread outages. One approach to mitigate impending cascading failures is through controlled islanding, in which a set of transmission lines is deliberately tripped to partition the unstable system into several disjoint, internally stable islands. Selecting such a set of transmission lines is inherently a combinatorial optimization problem. Current approaches address this problem in two steps: first classify coherent generators into groups and then separate generator groups into different islands with minimal load-generation imbalance. These methods, however, are based on computationally expensive heuristics that do not provide optimality guarantees. In this paper, we propose a novel approach to controlled islanding based on weak submodularity. Our formulation jointly captures the minimal generator non-coherency and minimal load-generation imbalance in one objective function. We relax the problem to a formulation with bounded submodularity ratio and a matroid constraint, and propose an approximation algorithm which achieves a provable optimality bound on non-coherency and load-generation imbalance. The proposed framework is tested on IEEE 39-bus and 118-bus power systems.

*Index Terms*—Cascading failure, controlled islanding, submodularity ratio, weak submodularity.

## I. Introduction

Due to the increasing demand for electricity and unpredictable supplies from renewable energy, power systems are being operated close to their stability limits. When a large disturbance such as a transmission line outage occurs in one geographic area, its nearby lines have to compensate for the failed component by carrying additional power flows [1]. This shift of power flows may cause other transmission lines to exceed their capacities, which leads to further outages. This process is referred to as *cascading failure*, in which local failures propagate and trigger successive failures throughout the entire power system [2]. Cascading failures can cause significant economic damage. For instance, in the 2003 North American blackout, two transmission line outages in the state of Ohio created a cascading failure that left 55 million people without power [3].

*Controlled islanding* is an effective corrective control to mitigate the impending cascading failure [4], [5]. Controlled islanding partitions the unstable power system into a set of smaller, disjoint subsystems, or islands, by tripping a selected set of transmission lines. Each of these islands is an internally stable and self-contained system, which can be reconnected to restore the entire system when the failure is cleared.

It is inevitable to shed a certain amount of load in some islands in order to match the limited generation within the island. A large amount of power imbalance between load and generation within an island may cause the power frequency to deviate from the acceptable region and hence lead to an eventual collapse of the island [5]. Thus, it is important to minimize the overall load-generation imbalance in each island. On the other hand, generators tend to oscillate at different frequencies during cascading failure. If the generators in an island do not oscillate at similar frequencies, the island may become unstable as well [6]. Thus, it is also critical to limit the range of generator frequencies within each island when choosing an islanding strategy.

Algorithms that select islanding (partition) strategies based on these considerations have been proposed in [5], [7], which address the problem in two steps. These methods first identify clusters of generators with similar rotor angle frequencies, which are referred to as coherent generator groups. After forming such groups, a set of transmission lines is selected to trip, such that the coherent groups of generators are separated into different islands while the load-generation imbalance within each island is minimized. Current two-step approaches to controlled islanding rely on graph-based heuristics such as min-cut and spectral clustering. Such approaches not only are unable to guarantee solutions' optimality but also lack the flexibility in making trade-off between generator coherency and load-generation imbalance, defined as the ability to sacrifice a small amount of generator coherency within each island in order to gain a large reduction in overall power imbalance.

In this paper, we propose a one-step approach to controlled islanding that jointly minimizes the generator non-coherency and power imbalance when selecting an islanding strategy, with provable optimality guarantees. Our key insight is that selecting a subset of transmission lines to trip in order to form desired partition is inherently a combinatorial problem with weak submodularity property. Weak submodularity is a weakened form of diminishing property of set functions, which is sufficient to guarantee the optimality of a greedy selection algorithm with a relaxed approximation factor [8]. We make the following specific contributions:

- We formulate the problem of selecting a set of transmission lines to trip in order to minimize the weighted sum of generator non-coherency and load-generation imbalance, subject to constraints that ensure that the reference generators are in disjoint islands.
- We relax the problem by removing the limits on load and generation and then prove that the relaxed formulation has

Z. Liu, L. Bushnell, D. Kirschen and R. Poovendran are with the Department of Electrical Engineering, University of Washington, Seattle, WA, 98195 USA (e-mail: {zhipliu, lb2, kirschen, rp3}@uw.edu).

A. Clark is with the Department of Electrical and Computer Engineering, Worcester Polytechnic Institute, Worcester, MA 01609 USA (e-mail: aclark@wpi.edu).

bounded submodularity ratio. We also prove a sufficient condition for the constraints using graph matroid.
- We propose a polynomial-time greedy algorithm to approximately solve the relaxed islanding problem and prove that the algorithm has optimality guarantees on the minimal non-coherency and load-generation imbalance.
- We propose a submodular maximization approach to selecting reference generators based on slow-coherency theory, which provides a optimality guarantee of $(1 - 1/e)$ on reference generator independence.
- We test our proposed approach to controlled islanding on the IEEE 39-bus and 118-bus systems. We validate the optimality of our method by comparing to a two-step spectral clustering based approach [7].

The rest of this paper is organized as follows. Section II reviews related work. Section III presents the power system model and background on submodularity and matroid. Section IV presents the problem formulation and the islanding algorithm. Section V discusses the slow-coherency based reference generator selection and proposes a submodular maximization approach to selecting reference generators. Section VI presents the simulation results. Section VII concludes the paper.

## II. Related Work

Controlled islanding has been studied extensively especially since the 2003 North American blackout that was caused by cascading failure [5], [7], [9], [10], [11], [12]. One important branch of study is slow coherency-based islanding, which was first introduced in [5]. This approach first classifies generators into groups with minimum non-coherency, and then finds a best cutset of transmission lines that separate generator groups into different islands while minimizing the overall load-generation imbalance through exhaustive search. This method, however, is not scalable to large systems. A refined islanding scheme is presented in [11], in which the step of searching for minimum-load-generation islands is approximated by finding minimal power-flow disruption between islands, and hence the problem is mapped to a classic graph min-cut problem. This allows the use of existing graph theory based tools, but an exhaustive search is still needed in order to find the optimal solution. In [9], the slow coherency-based generator grouping and min-cut approach to islanding have been implemented in a simulation study of the 2003 blackout scenario.

Other approaches to controlled islanding include *Ordered Binary Decision Diagram (OBDD)* methods [12] and machine learning based methods [10]. In [12], a large power system is first simplified into a topology with less than 30 nodes and then a fast brute-force search is used to find feasible cutsets. The topology simplification is heavily based on system operation heuristics, such as integrating important transmission lines that should not be cut. As a consequence, the OBDD approach cannot provide any guarantee on the optimality of solutions. Off-line machine learning based approaches have been investigated in [10]. Given a large set of off-line simulation data, a probability model can be established to assist evaluation of online islanding strategies. The optimality of such data-driven methods, however, is not guaranteed.

A two-step spectral clustering based control islanding is proposed in [7]. In the first step, the generators that have strong dynamic coupling are grouped, by applying the normalized spectral clustering to the system eigenbasis. In the second step, other buses are grouped into islands by constrained spectral clustering with the goal of minimizing power-flow disruption on cutting edges. Although this method is computationally efficient, the gap between the solution and the optimal strategy is not guaranteed.

Submodular optimization techniques have been established to address many power system stability issues, including voltage instability [13] and small signal instability [14]. A framework that provides optimality guarantees on the controlled islanding with scalable computational overhead, has not been studied in existing literature. A preliminary version of this work appeared in [15]. The optimality bound reported in [15] was derived based on a result from the existing submodularity literature [16] which was later shown to be incorrect in [17]. In this paper, we derive a new optimality bound via bounded submodularity ratio (weak submodularity) instead of submodularity. In addition to correcting these errors, we reformulate the two-step optimization approach to a one-step minimization problem, add a reference generator selection algorithm, and have a detailed validation of proposed islanding algorithms on both 39-bus and 118-bus power systems with comparison to the state-of-the-art spectral clustering method.

## III. System Model and Preliminaries

In this section, we present the power system model and introduce notations that will be used throughout the paper. We also give background on submodularity and matroids.

### A. Power System Dynamics and Generator Coherency

We consider a power system consisting of $n$ generators, $m$ buses and $l$ transmission lines that connect buses. Let $\delta_i$ denote the rotor angle of generator $i$ and $\Delta \delta_i$ denote the rotor angle deviation from a steady state operating point. Let $M_i$ denote the inertia of generator $i$. Given a base frequency $\omega_0$, the generator dynamics based on the classical linearized generator swing equation is given as [6]

$$\Delta \ddot{\delta} = M^{-1} K \, \Delta \delta, \tag{1}$$

where $\Delta \delta = [\Delta \delta_1, \ldots, \Delta \delta_n]^T$ is the state vector and $M = diag(2M_1/\omega_0, \ldots, 2M_n/\omega_0)$ is the inertia matrix. The matrix $K$ has entries

$$K_{ij} = \begin{cases} - V_i V_j B_{ij} \cos(\delta_i - \delta_j), & \text{if } i \neq j, \\ - \sum_{k=1, k\neq i}^{n} K_{ik}, & \text{if } i = j, \end{cases}$$

where $V_i$ is the per unit voltage behind transient reactance at generator $i$ and $B_{ij}$ is the imaginary part of the $(i,j)$-entry of the admittance matrix reduced to the internal generator nodes.

Following a large disturbance, generators often oscillate at different frequencies. Two synchronous generators are called $\epsilon$-coherent if the maximum change in their rotor angle difference falls within a specified tolerance, $\delta_{ij}(t) - \delta_{ij}(0) \leq \epsilon$, at any time $t$ [18]. We omit the $\epsilon$ for simplicity of notation. In the design of a strategy to mitigate impending cascading failures, it is critical to have coherent generators remain connected and separate any generators that are non-coherent into disjoint subsystems. Each subsystem is an island and the set of generators in the same island is called a generator group. Let $r$ be the desired number of islands.

Given a set of eigenvalues $\sigma_r = \{\lambda_1, \ldots, \lambda_r\}$ of the system matrix $M^{-1}K$, let $U$ be the eigenbasis of the $\sigma_r$ eigenspace. Two generators (states) are coherent with respect to $\sigma_r$ if their corresponding row vectors in $U$ are linearly dependent. With selection of the $r$ "slowest modes" of the system matrix $M^{-1}K$, representing the $r$ eigenvalues with smallest magnitudes, the generators dominant in each mode will have eigenbasis rows that are most linearly independent to each other [6]. Such a set of $r$ generators are called reference generators and can be used to initiate the $r$ islands in controlled islanding. The classic heuristics as well as our proposed submodular approach to reference generator selection are presented in Section V.

Given a set of reference generators, an $n \times r$ coherency matrix $L$ can be calculated, whose entries represent how coherent each generator is respect to a reference generator. The algorithm to obtain $L$ has the following steps [6]:

1) Choose the set of $r$ slowest modes, $\sigma_r$.
2) Compute a basis matrix $U_{n \times r}$ of the $\sigma_r$-eigenspace.
3) Let $U_1$ denote the matrix composed of the rows of $U$ corresponding to the reference generators.
4) Compute $L = U U_1^{-1}$.

Any generator grouping strategy can be represented by a partition matrix $L_g$ that has exactly one 1 per row and all other entries zero. Matrix $L_g$ has the same size as $L$. Each column of $L_g$ represents a generator group dominated by a reference generator. The $(i,j)$th entry of $L_g$ being 1 can be interpreted as that the generator $i$ is in the group of the $j$th reference generator. Given any grouping strategy, we define a metric that quantifies how coherent generators are within each group, given by

$$\overline{H}(L_g) = \|L - L_g\|_F^2, \qquad (2)$$

where $\|\cdot\|_F$ corresponds to the Frobenius norm. We note that $\overline{H}(L_g)$ measures the overall non-coherency of generators in each group and hence $\overline{H}(L_g) = 0$ when all generators within the same group are perfectly coherent ($\delta_{ij}(t) - \delta_{ij}(0) = 0, \forall t$).

In this subsection, we have established the dynamical model of the power system. The algebraic properties of the system will be modeled in the following subsection.

### B. Power flow network

Given a power system with $m$ buses and $l$ transmission lines, its topology can be depicted by a graph $G(N, E)$ where $N = \{1, \ldots, m\}$ is the index set of buses and $E \subseteq N \times N$ is the set of edges. We denote $e_k$ as the $k$th element in $E$. An edge $e_k = (i,j)$ is in $E$ if and only if there is a transmission line connecting bus $i$ and bus $j$. For any bus $i$, the set of its neighboring buses is defined by $\Delta(i) = \{j : (i,j) \in E\}$. A bus with a generator connected to it is called a generator bus while others are called load buses. Let $g_i$ be the active power generation and $d_i$ be the active power demand at bus $i$. The amount of active power flowing over a transmission line $(i,j) \in E$ from bus $i$ to $j$ is denoted by $p_{ij}$. At the initial operating point, we set $g_i = g_i^0$ and $d_i = d_i^0$ for all buses $i$ and $p_{ij} = p_{ij}^0$ for all transmission lines.

With transmission line resistance neglected, the net active power flowing into any bus equals the net active power demand at this bus, which is known as the flow conservation property:

$$\sum_{j \in \Delta(i)} p_{ji} = d_i - g_i, \ \forall i \in N. \qquad (3)$$

Using the incidence matrix $A$ of an orientation of the graph $G(N, E)$, the linear equations in (3) are equivalent to the matrix expression

$$Ax = b, \qquad (4)$$

where $A$ is an $m \times l$ matrix with each row corresponding to a bus and each column corresponding to a transmission line. The entry of $A$ in the $i$th row and $k$th column is denoted by $(A)_{ik}$, where $(A)_{ik} = -1$ if the line $e_k \in E$ connects bus $i$ to some bus $j$ with $i > j$, $(A)_{ik} = 1$ if $e_k \in E$ connects $i$ to $j$ with $i < j$, and $(A)_{ik} = 0$ otherwise. The vector $x \in \mathbb{R}^l$ consists of entries $x_k$, where the $k$th entry $x_k = p_{ij}$ corresponding to the power flow on the edge $e_k = (i,j) \in E$. The vector $b \in \mathbb{R}^m$ has entries $(d_i - g_i)$ corresponding to the net active power demand at each bus $i$.

Let $a_i$ be the $i$th column of the incidence matrix $A$ and for any subset $S \subseteq E$, define $A(S)$ to be a matrix consisting of columns $\{a_i, \forall e_i \in S\}$. Tripping any transmission line $e_k \in E$ will change the incidence matrix $A$ in Eq. (4) to $A(E \setminus \{e_k\}) = [a_1, \ldots, a_{k-1}, a_{k+1}, \ldots, a_l]$.

### C. Background on Submodularity and Matroids

Submodularity is a diminishing-returns property of set functions, wherein the incremental benefit of adding an element to a set $S$ decreases as more elements are added to $S$. Let $2^\Omega$ denote the power set of $\Omega$. A set function $f : 2^\Omega \to R$ is *submodular* if, for any sets $S \subseteq T \subseteq \Omega$ and any element $v \in \Omega \setminus T$, $f(S \cup \{v\}) - f(S) \geq f(T \cup \{v\}) - f(T)$.

The submodularity ratio measures to what extent a function $f$ has submodular properties, defined as follows.

**Definition 1** (Submodular Ratio, Weak Submodularity [19]). *Let $f : 2^\Omega \to \mathbb{R}$ be a non-negative set function. The submodularity ratio of $f$ with respect to a set $U \subseteq \Omega$ and a parameter $k \geq 1$ is*

$$\gamma_{U,k} = \min_{\substack{L, S : L \cap S = \emptyset \\ L \subseteq U, |S| \leq k}} \frac{\sum_{z \in S} (f(L \cup \{z\}) - f(L))}{f(L \cup S) - f(L)}$$

*$f$ is called weakly submodular at $U$ and $k$ if $\gamma_{U,k} > 0$.*

A matroid is defined as follows.

**Definition 2.** *Let $V$ denote a finite set, and let $\overline{\mathcal{I}}$ denote a collection of subsets of $V$. The tuple $\mathcal{M} = (V, \overline{\mathcal{I}})$ is a* matroid *if (i) $\emptyset \in \overline{\mathcal{I}}$, (ii) for any $B \in \overline{\mathcal{I}}$ and $A \subseteq B$, $A \in \overline{\mathcal{I}}$, and (iii)* *if $A, B \in \overline{\mathcal{I}}$ and $|A| < |B|$, then there exists $v \in B \setminus A$ such that $(A \cup \{v\}) \in \overline{\mathcal{I}}$.*

A maximal independent set of a matroid is a *basis*, denoted by $\mathcal{I}$. An important property of matroid bases is shown below.

**Lemma 1** ([20]). *If $\mathcal{I}_1$ and $\mathcal{I}_2$ are two bases of a matroid $\mathcal{M}(V, \overline{\mathcal{I}})$, then there exists a bijection $\pi : \mathcal{I}_1 \to \mathcal{I}_2$ such that*

$$\mathcal{I}_1 - x + \pi(x) \in \overline{\mathcal{I}} \text{ for all } x \in \mathcal{I}_1$$

An important class of matroids, denoted *graphic matroids*, is defined as follows.

**Lemma 2.** *Let $G = (N, E)$ be a graph with node set $N$ and edge set $E$. Define $\overline{\mathcal{I}}$ to be a collection of subsets of $E$ such that $A \in \overline{\mathcal{I}}$ iff $A$ does not contain any cycle (i.e., $A$ is a forest). Then $\mathcal{M} = (E, \overline{\mathcal{I}})$ is a matroid.*

For connected graphs, the set of bases of a graphic matroid is the set of all spanning trees.

## IV. Proposed Controlled Islanding Framework

This section formulates the problem of partitioning a power system into internally stable islands and presents an islanding algorithm with provable optimality guarantees by exploiting the weak submodularity of the problem.

### A. Problem Formulation

Suppose that we are given a power system with graph $G(N, E)$ and a set of reference generators at buses $\{s_1, \ldots, s_r\} \subseteq N$. We are also given a maximum generation limit $g_{\max}^i$ for each generator and a desired load $d_{\max}^i$ for each load at bus $i$. Note that $g_{\max}^i = 0$ for load buses and $d_{\max}^i = 0$ for buses without any load.

The problem studied in this paper is to select a set of edges $S' \subseteq E$ to trip such that the system is partitioned into $r$ disjoint internally stable islands where each island contains a unique reference generator. We denote $S$ as the complementary set of $S'$ such that $S = E \setminus S'$. Such a set $S$ is called a partitioning set. The partition of the system induced by $S$ is called an islanding strategy.

An island being internally stable requires that the generators within the island must be coherent and the imbalance between load and generation at all buses must be minimized [5]. The load-generation imbalance is an estimate of how much load must be shed during controlled islanding. In the existing literature [7], [11], the load-generation imbalance is often defined as sum of mismatch between load and generation within each island or sum of power flows on lines to be tripped, based on a power flow solution at a steady state prior to the disturbance. For islands having power imbalance, the load that must be shed depends on the new steady state established in the island.

A metric that captures the minimum load that must be shed, given any partitioning set $S$, is defined as follows

$$F(S) \triangleq \begin{aligned} &\min_{x,d,g} && \|A(S)x - (d_0 - g_0)\|_2^2 \\ &\text{s.t.} && 0 \leq d \leq d_{\max}, \ 0 \leq g \leq g_{\max}, \\ & && A(S)x = d - g, \end{aligned} \quad (5)$$

where $A(S)$ is the incidence matrix induced by the set of edges $S$. The vectors $d_{\max} = [d_{\max}^1, \ldots, d_{\max}^m]^T$ and $g_{\max} = [g_{\max}^1, \ldots, g_{\max}^m]^T$ are the maximum available load and generation at each bus, respectively. The vector $(d_0 - g_0) = [d_1^0 - g_1^0, \ldots, d_m^0 - g_m^0]^T$ is the net load at initial operating point. The metric $F(S)$ is equal to the minimum load-generation imbalance for any power flow $x$ that satisfies the flow conservation constraint (4) at each bus.

Generator coherency is another essential property to achieve the internal stability of islands. In the following, we define a metric based on (2) that evaluates the generator non-coherency given any partitioning set $S$. By linearizing the system dynamics at a steady-state operating point, a model can be established as described by (1). Following the algorithm in Section III-A, one can obtain a coherency matrix $L$ with respect to the given reference generators. Given any set $S$, a partition matrix $L_g$ can be calculated and the overall generator non-coherency is measured by $\overline{H}(L_g)$ in Eq. (2).

Let $u_i$ denote the bus index of the $i$th generator. For each non-reference generator at bus $u_i \in N \setminus \{s_1, \ldots, s_r\}$, define a function $H_i$ to be the optimal value of a quadratic program:

$$H_i(S) \triangleq \begin{aligned} &\min_{x} && \|A(S)x - c^i\|_2^2 \\ &\text{s.t.} && (A(S)x)_j = 0, \ \forall j \notin \{u_i, s_1, \ldots, s_r\}, \end{aligned} \quad (6)$$

where $A(S)$ is the incidence matrix induced by $S$ and $c^i$ is a vector defined by

$$(c^i)_j = \begin{cases} 1, & j = u_i \\ -L_{ik}, & j = s_k \in \{s_1, \ldots, s_r\} \\ 0, & \text{else.} \end{cases}$$

The following lemma establishes the relation between the metrics $H_i$ and $\overline{H}(L_g)$.

**Lemma 3.** *If the set of edges $S$ induces a partition such that each reference generator $s_k \in \{s_1, \ldots, s_r\}$ is in a unique island and the generator at bus $u_i \in N \setminus \{s_1, \ldots, s_r\}$ is connected to exactly one reference generator at bus $s_j$, then*

$$H_i(S) = \frac{1}{2}(1 - L_{ij})^2 + \frac{1}{2}\sum_{k=1, k\neq j}^{r} L_{ik}^2 = \frac{1}{2}\|(L - L_g)_i\|_2^2,$$

*where $(L - L_g)_i$ is the $i$th row of the matrix $(L - L_g)$.*

*Proof.* Suppose that the set of edges $S$ induces a partition of the system $G(N, E)$, where each reference generator $s_k \in \{s_1, \ldots, s_r\}$ is in a unique island and the generator at bus $u_i \in N \setminus \{s_1, \ldots, s_r\}$ is connected to a reference generator $s_j$. With an orientation of the graph $G$, the matrix $A(S)$ has the following block pattern

$$A(S) = \begin{bmatrix} A_1 & & \\ & \ddots & \\ & & A_r \end{bmatrix},$$

where each block $A_k$ corresponds to the incidence matrix of an island containing reference generator $s_k$. Similarly, the vectors in (6) become $x = [\mathbf{x}_1, \ldots, \mathbf{x}_r]^T$ and $c^i = [c_1^i, \ldots, c_r^i]^T$. Thus,

$$H_i(S) = \sum_{k=1}^{r} \min_{\mathbf{x}_k} \|A_k \mathbf{x}_k - c_k^i\|_2^2, \quad (7)$$

subject to $(A(S)x)_k = 0, \forall k \notin \{u_i, s_1, \ldots, s_r\}$.

Each term of $H_i(S)$ in (7) is a minimization problem. The solution to each problem is equivalent to the load-generation imbalance defined by (5) of a unique island. For the island $k = j$ which contains the bus $u_i$, given a "generation" of 1 at $u_i$ and a "load" of $L_{ij}$ at $s_j$, it has a "load-generation imbalance" of $(1 - L_{ij})^2/2$ resulting from the corresponding minimization problem. For each of the rest islands where $k \neq j$, there is only a "load" of $L_k$ at bus $s_k$ without any generation in the island, which results in a "load-generation imbalance" of $L_{ik}^2$. With the analogy above, we get

$$H_i(S) = \frac{1}{2}(1 - L_{ij})^2 + \frac{1}{2} \sum_{k=1, k \neq j}^{r} L_{ik}^2.$$

$L_g$ is the generator partitioning matrix induced by $S$. If the generator at bus $u_i$ is in the same island of the reference generator $s_j$, then the $i$th row of $L_g$ has a 1 at $j$th column and 0 anywhere else. Thus, we have

$$\|(L - L_g)_i\|_2^2 = (1 - L_{ij})^2 + \sum_{k=1, k \neq j}^{r} L_{ik}^2,$$

and hence $H_i(S) = \|(L - L_g)_i\|_2^2/2$, completing the proof. □

The proof of Lemma 3 is established using the analogy of the formulation in (5), in which we have a "generation" of 1 at bus $u_i$ and a "load demand'" of $L_{ik}$ at bus $s_k$, $\forall k \in \{1, \ldots, r\}$. We restrict $d_{\max}$ and $g_{\max}$ to be 0 for buses not in the set $\{u_i, s_1, \ldots, s_r\}$.

The function $H_i(S)$ can be interpreted as the sum of non-coherency in each group caused by the $i$th generator, given an islanding strategy $S$. We note that $H_i(S) = 0$ if the generator $i$ is a reference generator.

By Lemma 3 and the fact that $\overline{H}(L_g) = \|L - L_g\|_F^2 = \sum_{i=1}^{n} \|(L - L_g)_i\|_2^2$, we find that $\overline{H}(L_g) = 2\sum_{i=1}^{n} H_i(S)$.

The problem is to select a set of edges, $S \subseteq E$, to partition the system $G(N, E)$ into $r$ disjoint islands, $G(N_1, S_1), \ldots, G(N_r, S_r)$. Each given reference generator $s_k \in \{s_1, \ldots, s_r\}$ must be in a unique island, $s_k \in N_k, \forall k$, while minimizing the overall load-generation imbalance and generator non-coherency, which can be formulated as

$$\begin{aligned} \min_{S} \quad & \xi F(S) + \sum_{i=1}^{n} H_i(S) \\ \text{s.t.} \quad & S = S_1 \cup \cdots \cup S_r, N = N_1 \cup \cdots \cup N_r, \\ & s_k \in N_k, \ \forall k \in 1, \ldots, r \end{aligned} \quad (8)$$

with a trade-off parameter $\xi \geq 0$ balancing the two objectives.

In the following section, we exploit the weak submodularity of problem (8), which leads to polynomial-time approximation algorithms with provable bound on the optimality of solutions.

*B. Weak Submodularity Approach*

In order to define a relaxation of the problem (8) which has weak submodularity, we introduce the following metrics

$$h_i(S) = \min_x \|A(S)x - c^i\|_2^2,$$

for all $i \in \{1, \ldots, n\}$ and

$$f(S) = \min_x \|A(S)x - b_0\|_2^2,$$

where $b_0 = d_0 - g_0$ and the restrictions on variables $d$, $g$ and $x$ are removed. The following lemma establishes the connection between metrics $f(S)$ and $F(S)$, as well as $h_i(S)$ and $H_i(S)$.

**Lemma 4.** *Given any set $S \subseteq E$, we have $f(S) = 0$ if $F(S) = 0$, and $f(S) \leq F(S)$ if $F(S) > 0$; we have $h_i(S) = 0$ if $H_i(S) = 0$, and $h_i(S) \leq H_i(S)$ if $H_i(S) > 0$, for all $i$.*

*Proof.* If $F(S) = 0$, then there exists a power flow $\bar{x}$ such that $A(S)\bar{x} = b_0$. Thus, we have

$$f(S) = \min_x \|A(S)x - b_0\|_2^2 = \|A(S)\bar{x} - b_0\|_2^2 = 0.$$

If $F(S) > 0$, let $\hat{x}$ be the optimal solution of $F(S)$, i.e., $\hat{x} = \arg\min_x \|A(S)x - (d_0 - g_0)\|_2^2$ s.t. $0 \leq d \leq d_0, 0 \leq g \leq g_0$ and $A(S)x = d - g$. Then we have

$$f(S) = \min_x \|A(S)x - b_0\|_2^2 \leq \|A(S)\hat{x} - b_0\|_2^2 = F(S).$$

Similar proofs for $h_i(S)$ and $H_i(S)$ can be established following the same steps. □

By Lemma 4, the function $f(S)$ is a lower bound to the load-generation imbalance $F(S)$ by relaxing the restrictions on generation and load at each bus, while $h_i(S)$ is a lower bound to the generator non-coherency $H_i(S)$ for each generator $i$.

Thus, the problem (8) has a relaxed formulation, given as follows, that provides a lower bound on the optimal solution:

$$\begin{aligned} \min_{S} \quad & \xi f(S) + \sum_{i=1}^{n} h_i(S) \\ \text{s.t.} \quad & S = S_1 \cup \cdots \cup S_r, N = N_1 \cup \cdots \cup N_r, \\ & s_k \in N_k, \ \forall k \in 1, \ldots, r \end{aligned} \quad (9)$$

The following theorem establishes the weak submodularity structure of the objective function in (9).

**Theorem 1.** *Let $J(S) = \xi f(S) + \sum_{i=1}^{n} h_i(S)$. Denote $\gamma_{U,k}$ as the submodularity ratio of $J(S)$. Then for any set $U \subseteq E$ and $k \geq 1$, the submodularity ratio $\gamma_{U,k}$ is bounded by*

$$\gamma_{U,k} \geq \lambda_{\min}(C, k + |U|) \geq \lambda_{\min}(C),$$

*where $C = A^T A/(2n)$; $\lambda_{\min}(C)$ is the smallest eigenvalue of $C$ and $\lambda_{\min}(C, k)$ is the smallest $k$-sparse eigenvalue of $C$, defined by $\lambda_{\min}(C, k) = \min_{S:|S|=k} \lambda_{\min}(C(S))$ where $C(S) = A(S)^T A(S)/(2n)$.*

*Proof.* For any vector $v$, we define a metric

$$g_v(S) = \min_x \|A(S)x - v\|_2^2.$$

Equivalently, $g_v(S) = \min_y \|(A(S)/\sqrt{2}) y - v\|_2^2$ where each column of the matrix $(A(S)/\sqrt{2})$ has norm 1 and sum of entries 0. The metric $g_v(S)$ has been shown to have bounded submodularity ratio $\gamma'_{U,k} \geq \lambda_{\min}(C, k + |U|) \geq \lambda_{\min}(C)$ [8], [19], [21]. Hence, it suffices to show that $J(S)$ has submodularity ratio $\gamma_{U,k}$ bounded by $\gamma'_{U,k}$, i.e., $\gamma_{U,k} \geq \gamma'_{U,k}$.

Let $\{v_i\}$ be a set of vectors where $v_0 = b_0$ and $v_i = c^i$ for $i \in \{1, \ldots, n\}$. Let $\{a_i\}$ be a set of weights where $a_0 = \xi$ and $a_i = 1$ for $i \in \{1, \ldots, n\}$. Then $J(S) = \sum_{i=0}^{n} a_i g_{v_i}(S)$.

$$\gamma_{U,k} = \min_{\substack{L, S: L \cap S = \emptyset \\ L \subseteq U, |S| \le k}} \frac{\sum_{z \in S} (J(L \cup \{z\}) - J(L))}{J(L \cup S) - J(L)}$$

$$= \min_{\substack{L, S: L \cap S = \emptyset \\ L \subseteq U, |S| \le k}} \frac{\sum_i \sum_{z \in S} (a_i g_{v_i}(L \cup \{z\}) - a_i g_{v_i}(L))}{\sum_i (a_i g_{v_i}(L \cup S) - a_i g_{v_i}(L))}$$

$$\ge \min_{\substack{L, S: L \cap S = \emptyset \\ L \subseteq U, |S| \le k}} \left( \min_i \frac{\sum_{z \in S} (g_{v_i}(L \cup \{z\}) - g_{v_i}(L))}{(g_{v_i}(L \cup S) - g_{v_i}(L))} \right)$$

$$\ge \gamma'_{U,k}$$

completing the proof. □

In what follows, we define a graphic matroid constraint that is sufficient to satisfy the constraint on $S$ of problem (9). Given the graph $G(N, E)$ and reference generators at buses $\{s_1, \ldots, s_r\} \in N$, we construct a new graph $G_0(N_0, E_0)$ as follows. The vertex set $N_0 = N \cup \{0\}$ and the edge set $E_0 = E \cup \{(0, s_i), i = 1, \ldots, r\}$. Intuitively, the graph $G_0$ is obtained by adding a new node 0 to the graph $G$ and creating new edges to connect all reference generators to node 0.

Let $\mathcal{M} = (E_0, \mathcal{I}_0)$ denote the graphic matroid defined on the graph $G_0$, with basis $\mathcal{I}$ being the collection of all spanning trees of $G_0$. The following lemma relates the basis $\mathcal{I}$ of $\mathcal{M}$ to the constraint in (9).

**Lemma 5.** *Given a graph $G(N, E)$ and a set of reference generators $\{s_1, \ldots, s_r\} \in N$, a set of edges $S \in E$ induces a partition of the graph: $S = \{\cup_{k=1}^{r} S_k\}$, $N = \{\cup_{k=1}^{r} N_k\}$ and $S_i \cap S_j = \emptyset, N_i \cap N_j = \emptyset, \forall i \ne j$, where each reference generator satisfies $s_k \in N_k, \forall k$, if*

$$S \cup \{(0, s_i), i = 1, \ldots, r\} \in \mathcal{I}. \quad (10)$$

*Proof.* The proof is established by contradiction. If $S \cup \{(0, s_i), i = 1, \ldots, r\} \in \mathcal{I}$, then the set of edges $S \cup \{(0, s_i), \forall i\}$ forms a spanning tree of $G_0$. The node 0 is the root of the tree and $S = \{\cup_{k=1}^{r} S_k\}$ and $N = \{\cup_{k=1}^{r} N_k\}$ form $r$ branches of the tree. Given any two reference generators at nodes $s_i$ and $s_j$, then they are connected by a path $(s_i, 0), (0, s_j)$ in graph $G_0$. Suppose that two reference generators are in the same branch, i.e., $s_i \in N_i$ and $s_j \in N_i$. Then $S$ contains edges that connect $s_i \in N_i$ and $s_j \in N_i$ by a path $(s_i, v_1), (v_1, v_2), \ldots, (v_k, s_j) \subseteq S$ in graph $G$. Since node $0 \notin N$, we have $0 \ne v_1 \ne \cdots \ne v_k$. Thus, the path $(0, s_i), (s_i, v_1), (v_1, v_2), \ldots, (v_k, s_j), (s_j, 0) \subseteq S \cup \{(0, s_i), \forall i\}$ forms a cycle in graph $G_0$, which contradicts with the definition of a spanning tree. Therefore, by contradiction, we have that $s_k \in N_k, \forall k$. □

By Lemma 5, the problem of selecting a set of edges $S \subseteq E$ to form disjoint islands where each island contains a unique reference generator, while minimizing load-generation imbalance and generator non-coherency can be formulated as

$$\begin{aligned} \min_S \quad & J(S) \\ \text{s.t.} \quad & S \cup \{(0, s_i), i = 1, \ldots, r\} \in \mathcal{I}. \end{aligned} \quad (11)$$

where $J(S) = \xi f(S) + \sum_{i=1}^{n} h_i(S)$.

**Algorithm 1** Algorithm for selecting a set of edges that partitions the system into internally stable islands.

1: **procedure** CONTROL_ISLAND($G(N, E), \{s_i\}, b_0, L, \xi$)
2:  **Input**: Power system topology $G(N, E)$, set of reference generators $\{s_i\} \subseteq N$, net load at each bus $b_0$, generator coherency matrix $L$, trade-off parameter $\xi$.
3:  **Output**: Set of edges $S \in E$ that form a partition of $G$, where each partition contains a group of coherent generators while minimizing load-generation imbalance.
4:  **Initialization:** Construct metric $J = \xi f + \sum_{i=1}^{n} h_i$; Construct new graph $G_0(N_0, E_0)$ by adding new node 0 and edges $(0, s_i), \forall i$ to $G$; Set $\Omega = E, S = \emptyset$.
5:  *1) Greedy Selection:*
6:  **while** $\Omega \ne \emptyset$ **do**
7:   $e \leftarrow \arg\max_{v \in \Omega \setminus S} (J(S) - J(S \cup \{v\}))$
8:   **if** $S \cup \{e\} \cup \{(0, s_i), \forall i\}$ contains no cycle **then**
9:    $S \leftarrow S \cup \{e\}$
10:   **end if**
11:   $\Omega \leftarrow \Omega \setminus \{e\}$
12:  **end while**
13:  *2) Local Search:*
14:  **while** $\exists v \in S, \exists e \in E \setminus S : J(S \setminus \{v\} \cup \{e\}) < (1-\epsilon) J(S)$ and $S \setminus \{v\} \cup \{e\} \cup \{(0, s_i), \forall i\}$ has no cycle **do**
15:   $S \leftarrow S \setminus \{v\} \cup \{e\}$
16:  **end while**
17:  **return** $S$
18: **end procedure**

### C. Proposed Islanding Algorithm

In this subsection, we present a greedy algorithm that approximately solves (11) with provable optimality guarantees. The algorithm consists of two stages, namely greedy selection and local search. Each stage proceeds as follows. In greedy selection, the set $S$ is initialized to be empty and the groundset $\Omega$ is initialized to be the set of edges $E$. At each iteration, the algorithm selects an edge $e \in \Omega \setminus S$ that maximizes $J(S) - J(S \cup \{e\})$ and deletes $e$ from $\Omega$. The edge $e$ is added to the set $S$ if it does not generate any cycle in the graph consisting of edges $S \cup \{e\} \cup \{(0, s_i), \forall i\}$. The greedy selection terminates when the groundset $\Omega$ is empty.

In local search, let $\epsilon > 0$ be a constant parameter. The algorithm selects a pair of edges $(v \in S, e \in E \setminus S)$ that by swapping them can lower the function $J$ value, i.e., $J(S \setminus \{v\} \cup \{e\}) < (1-\epsilon) J(S)$, and the graph consisting of edges $S \setminus \{v\} \cup \{e\} \cup \{(0, s_i), \forall i\}$ contains no cycle. A pseudocode description is given as Algorithm 1.

Intuitively, the algorithm first selects a spanning tree from graph $G_0$ while minimizing the function $J$. Then in local search, the algorithm swaps edges from the tree to converge to a local optimum which further minimizes $J$. The following proposition describes the optimality guarantee of Algorithm 1.

**Proposition 1.** *Let $S^*$ denote the optimal solution to (11) and let $S$ be the solution returned by Algorithm 1. Then, we have*

$$J(S) \le (m - r - \gamma_0) J(S_{t-1}) + \gamma_0 J(S^*),$$

where $\gamma_0 = \lambda_{\min}(C, 2|S|)$ and $S_{t-1}$ denotes the selected set at the second-to-last iteration of Algorithm 1.

*Proof.* The matroid $\mathcal{M}$ has rank $m$ and hence $|S^*| = |S| = m - r$. Let $e_i$ be the $i$th element selected by the Algorithm 1 and let $S_i = \{e_1, \ldots, e_i\}$ be the set containing the first $i$ elements picked by the algorithm. Let $\pi : S^* \to S$ be a bijection as in Lemma 1 so that the elements $e_i^*$ of $S^*$ are sorted to satisfy $S \setminus \{e_i\} \cup \{e_i^*\} \in \overline{\mathcal{I}}$ for all $i$.

By the definition of submodularity ratio, we have

$$\frac{\sum_{z \in S^*} (J(S_{i-1}) - J(S_{i-1} \cup \{z\}))}{J(S_{i-1}) - J(S_{i-1} \cup S^*)} \geq \gamma_{S,|S|} \geq \gamma_0, \ \forall i,$$

where $\gamma_0 = \lambda_{\min}(C, 2|S|)$. By the greediness of Algorithm 1 and that $S_{i-1} \cup \{e_i^*\} \in \overline{\mathcal{I}}$ for all $i$, we have

$$J(S_{i-1}) - J(S_i) + \sum_{z \in S^* \setminus e_i^*} (J(S_{i-1}) - J(S_{i-1} \cup \{z\}))$$
$$\geq \gamma_0 (J(S_{i-1}) - J(S_{i-1} \cup S^*)).$$

Since $J$ is monotone decreasing, we have

$$J(S_{i-1}) - J(S_i) + \sum_{z \in S^* \setminus e_i^*} J(S_{i-1}) \geq \gamma_0 (J(S_{i-1}) - J(S^*)),$$

or equivalently,

$$|S^*| J(S_{i-1}) - J(S_i) \geq \gamma_0 (J(S_{i-1}) - J(S^*)),$$

and hence

$$J(S_i) \leq (|S^*| - \gamma_0) J(S_{i-1}) + \gamma_0 J(S^*).$$

At the last iteration, we have

$$J(S) \leq (|S^*| - \gamma_0) J(S_{t-1}) + \gamma_0 J(S^*),$$

where $S_{t-1}$ denotes the selected set at the second-to-last iteration of Algorithm 1. □

The complexity of Algorithm 1 is described by the following proposition.

**Proposition 2.** *The runtime of Algorithm 1 is bounded above by $O(l^5 n + l^2 C)$, where the greedy selection stage requires a computation of at most $O(l^5 n)$ and the local search requires $O(l^2 C)$, where*

$$C = \frac{\log J(E) - \log J(\emptyset)}{\log(1 - \epsilon)}.$$

*Proof.* The runtime of Algorithm 1 is determined by the complexity of evaluating the objective function $J(S)$ and checking the graphical matroid constraint. For greedy selection stage, the algorithm requires $O(l^2)$ evaluations of $J(S)$, each of which involves solving $(n + 1)$ least-squares problems with complexity of $O(l^3)$, and $O(l)$ checks of the matroid constraints. Checking the matroid constraint is to detect cycle in an undirected graph which requires computation of $O(l)$. Thus, the total complexity for greedy selection is $O(l^5 n)$.

In local search, let $S_i$ denote the set $S$ after $i$ iterations of the algorithm. By construction, $J(S_i) < (1 - \epsilon) J(S_{i-1})$ and hence $J(S_i) < (1-\epsilon)^i J(S_0)$. Let $k$ be the number of iterations before the algorithm terminates. By the monotonicity of $J$, we have $J(E) = \min_S J(S)$ and $J(\emptyset) = \max_S J(S)$. Therefore,

$$J(E) \leq J(S_k) < (1 - \epsilon)^k J(S_0) \leq (1 - \epsilon)^k J(\emptyset).$$

Rearranging terms yields $k \leq C$. Each iteration requires at most $l^2$ evaluations of the objective function, implying that the worst-case computation is $O(l^2 C)$. □

## V. Reference Generator Selection

In this section, we discuss the slow-coherency based reference generator selection and present a submodular optimization approach to selecting reference generators.

From the discussion in Section III-A, given the set of $r$ slowest modes $\sigma_r$ and its corresponding eigenbasis matrix $U$, two generators (states) are coherent with respect to $\sigma_r$ if their corresponding row vectors in $U$ are linearly dependent. Then selecting $r$ reference generators is equivalent to find $r$ rows from $U$ that are most linearly independent.

We define a matrix $U(T)$ consisting of the selected set $T$ of rows from $U$. A measure of the linear independence of the row vectors of $U(T)$ is given by the Gramian

$$G_U(T) = \det(U(T)U(T)'),$$

where $U(T)'$ is the transpose of $U(T)$.

A physical interpretation of the Gramian $G_U(T)$ is that the row vectors of $U(T)$ forms a parallelepiped in the case of $|S| = 3$ and $G_U(T)$ is the volume of the parallelepiped. For $|T| > 3$, the parallelepiped becomes a $|T|$-dimensional polytope formed by the row vectors of $U(T)$.

The linear independence of the rows of $U(T)$ reaches maximum when the volume of the parallelepiped/polytope formed by $U(T)$ is largest. Therefore, the reference generator selection can be formulated as selecting a set $T$ of rows from $U$ such that

$$\max_{T \subseteq \{1, \ldots, n\}} G_U(T), \quad s.t. \ |T| = m. \quad (12)$$

A classic heuristic method to approximately solve (12) is by applying complete-pivoting Gaussian elimination to matrix $U$ [6]. Following the elimination, the generator (state) dominant in each column is chosen as a reference generator. This heuristic, however, does not have any guarantee on the optimality of solutions.

In what follows, we present a submodular optimization based approach to reference generator selection, which provides a provable guarantee on the maximum linear independence achieved.

Since $U(T)U(T)'$ is positive semi-definite, its determinant $G_U(T)$ is nonnegative. By the fact that $\log(\cdot)$ is a monotone increasing function, the optimal point of $\max_T G_U(T)$ is the same as $\max_T \log G_U(T)$. Thus, the problem (12) is equivalent to

$$\max_T \ \log \det(U(T)U(T)'), \quad s.t. \ |T| = m. \quad (13)$$

**Theorem 2.** *The function $\phi(T) = \log \det(U(T)U(T)')$ is submodular in $T$.*

**Algorithm 2** Algorithm for selecting reference generators.

1: **procedure** REFERENCE_GEN($\hat{A}, r$)
2:    **Input:** Matrix $\hat{A} = M^{-1}K$, number of islands $r$
3:    **Output:** Set of reference generators $T$
4:    **Initialization:** Calculate the set of $r$ slowest modes $\sigma_r$ and basis matrix $U$; construct the set function $\phi$; $T \leftarrow \emptyset$
5:    **while** $|T| < r$ **do**
6:      **for** $v \in \Omega \setminus T$ **do**
7:        $\delta_v \leftarrow \phi(T \cup \{v\}) - \phi(T)$
8:      **end for**
9:      $v^* \leftarrow \arg\max_v \delta_v$
10:     $T \leftarrow T \cup \{v^*\}$
11:   **end while**
12:   **return** $T$
13: **end procedure**

*Proof.* Let $\Sigma = UU'$. Suppose $X \in \mathbb{R}^r$ is multivariate Gaussian with probability density function

$$p(x) = \frac{1}{\sqrt{|2\pi\Sigma|}} \exp\left(-\frac{1}{2}(x-u)'\Sigma^{-1}(x-u)\right)$$

with some mean vector $u$. The (differential) entropy of the random variable $X$ is given by

$$\Phi(X) = \log\sqrt{|2\pi\Sigma|} = \log\sqrt{(2\pi)^r|\Sigma|}.$$

For a subset $T$ of variables, $X_T$, we denote $\Sigma_T = V(T)V(T)'$ and we have

$$\Phi(X_T) = \log\sqrt{(2\pi)^{|T|}|\Sigma_T|} = \frac{1}{2}|T|\log 2\pi + \frac{1}{2}\log|\Sigma_T|,$$

where the term $(|T|\log 2\pi)$ is modular in $T$.

Since entropy $\Phi(X_T)$ is submodular in $T$, the term $\log|\Sigma_T| = \log\det(V(T)V(T)') = \phi(T)$ is submodular. □

By Theorem 2, the problem (13) is submodular maximization subject to a cardinality constraint. A simple greedy algorithm can approximately solve (13) in polynomial time $O(nr)$ up to a worst-case optimality bound of $(1-1/e)$ [22]. We present such an algorithm, described by Algorithm 2.

## VI. NUMERICAL STUDY

This section presents the simulation results of our proposed approach to controlled islanding on IEEE 39-bus system and IEEE 118-bus system, with comparison to the state-of-the-art two-step spectral clustering based islanding method [7].

### A. System Setup

For both test cases, the initial system data including topology, transmission line impedance, load and generator configuration are specified by the test system. We obtain an initial operating point including knowledge of generator voltage magnitudes $V_i$ and angles $\delta_i$, by solving an optimal power flow using *Matpower* [23]. The generator inertia $M_i$ for the 39-bus and the 118-bus systems are given by [24] and [25], respectively. We set the desired number of islands $r = 3$.

We evaluate the performance of an islanding strategy $S$ by calculating the resulting load-generation imbalance and generator non-coherency. We adopt the metric $\sqrt{f(S)} = \text{dist}(b_0, \text{span}(A(S)))$ for estimating the load-generation imbalance with unit MW and $\overline{H}(L_g) = \|L - L_g\|_F^2$ for generator non-coherency where $L_g$ is the generator partition matrix induced by $S$. The trade-off parameter $\xi$ is chosen from the interval $[0, 1]$, in order to take the value of load-generation imbalance down to the same level of generator non-coherency.

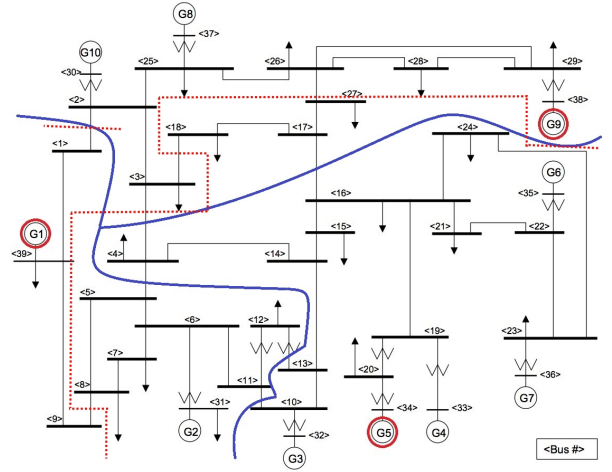

Fig. 1: IEEE 39-bus system topology. The circled generators {G1,G5,G9} are reference generators identified by Algorithm 2. The solid curves represent the cutset selected by our Algorithm 1 with parameter $\xi = 10^{-7}$. The dotted lines represent the cutset selected by the spectral clustering method.

The spectral clustering based approach to controlled islanding is implemented in two steps. Denote $\Omega_G \subseteq N$ as the set of generator buses. The first step is to classify generators into two coherent groups $T_1, T_2 \subseteq \Omega_G$ by minimizing the dynamic coupling between groups:

$$\min_{T_1, T_2 \subseteq \Omega_G} \sum_{j \in T_2} \sum_{i \in T_1} \left(\frac{\partial p_{ij}}{\partial \delta_{ij}}\left(\frac{1}{M_i} + \frac{1}{M_j}\right)\right), \quad (14)$$

where $\partial p_{ij}/\partial \delta_{ij} = |V_iV_jB_{ij}\cos(\delta_i - \delta_j)|$. In the second step, the algorithm finds a set of edges to cut such that the system is split into two subsystems $S_1, S_2 \subseteq N$ where $T_1 \subseteq S_1$ and $T_2 \subseteq S_2$, while minimizing the power-flow disruption on cutting edges:

$$\min_{S_1, S_2 \subseteq N} \sum_{i \in S_1, j \in S_2} |p_{ij}|, \quad s.t. \; T_1 \subseteq S_1, T_2 \subseteq S2. \quad (15)$$

For $r > 2$, repeat (14) and (15) for each subsystem until the desired number of islands $r$ is achieved.

### B. Test Case I : IEEE 39-bus Test System

The IEEE 39-bus test system contains 10 generators and 46 transmission lines. The system topology is shown in Fig. 1. Algorithm 2 identifies generators {G1, G5, G9} as reference generators, which is consistent with the result of the classic Gaussian elimination based reference generator selection.

By choosing $\xi = 10^{-7}$, Algorithm 1 finds the optimal islanding strategy based on metric $J(S)$ by tripping lines {1-2,

TABLE I: Comparison of our proposed islanding ($\xi = 10^{-7}$) and spectral clustering islanding results on the 39-bus system.

| Method | $J(S)$ | Load-generation Imbalance $\sqrt{f(S)}$ | Generator Non-coherency $\overline{H}$ |
|---|---|---|---|
| Proposed Islanding | 0.1696 | 241.1 MW | 1.4237 |
| Spectral Clustering | 0.2149 | 330.3 MW | 1.3104 |

TABLE II: Comparing results of the proposed islanding methods with different $\xi$ values on the 39-bus system.

| Trade-off Parameter | Islands[1] | Imbalance $\sqrt{f(S)}$ | Non-coherency $\overline{H}$ |
|---|---|---|---|
| $\xi = 0$ | {1,4,5,7,8,9,39}<br>{2,3,17,18,25,26,27,28,29,30,37,38}<br>{6,10:16,19:24,31:36} | 758.3 MW | 1.3104 |
| $\xi > 10^{-5}$ | {2,3,4,5,25,26,27,28,29,30,37,38}<br>{1,6,7,8,9,10,11,31,32,39}<br>{12:24,33,34,35,36} | 28.5 MW | 1.5900 |

3-4, 4-5, 10-11, 12-13, 16-17}. The coherent generator groups found are {G1,G2}, { G3,G4,G5,G6,G7}, {G8,G9,G10}.

The two-step spectral clustering method finds generator groups {G1}, {G2,G3,G4,G5,G6,G7}, {G8,G9,G10} in the first step and a min-cut set of edges {1-2, 8-9, 3-4, 3-18, 17-27} in the second step.

Fig. 1 illustrates the difference between islanding strategies of our proposed approach and the spectral clustering method on topology. A comparison of performance between the two islanding strategies is shown in Table I. We observe that our strategy reduces the load-generation imbalance by 27% compared to the spectral clustering islanding, although the generator non-coherency is 8.6% higher.

Table II shows islanding results from Algorithm 1 with different $\xi$ values. When choosing $\xi = 0$, only generator non-coherency is considered in the problem (11) and hence the objective function $J(S) = \sum_{i=1}^{n} h_i(S)$. The resulting submodular islanding strategy is able to find the same coherent generator groups as the spectral clustering method. When $\xi > 10^{-5}$, we have $J(S) \approx \sum_{i=1}^{n} h_i(S)$. The objective function $J(S)$ is dominated by the load-generation imbalance, while the effect of generator non-coherency is negligible. The Algorithm 1 finds a near-balanced islanding strategy by sacrificing 21.3% on generator non-coherency compared to the spectral clustering method.

*C. Test Case II : IEEE 118-bus Test System*

For this test case, the system topology is given by Fig. 2. Generators at buses {10, 54, 87} are selected by Algorithm 2 as reference generators, which is the same as classic reference generator selection result.

When $\xi = 10^{-6}$, Algorithm 1 returns a solution that has the same generator non-coherency but a lower load-generation imbalance compared to the spectral clustering solution, as shown in Table III. A comparison between the two solutions on system topology is shown in Fig. 2.

We conclude that our proposed islanding approach can identify solutions that achieve lower load-generation imbalance compared to the two-step spectral clustering method, with the same or slightly higher generator non-coherency. This is due to the trade-off mechanism of the proposed approach that is able to reduce the load-generation imbalance by changing coherent generator groups.

---
[1] {10:16} represents {10,11,12,13,14,15,16}.

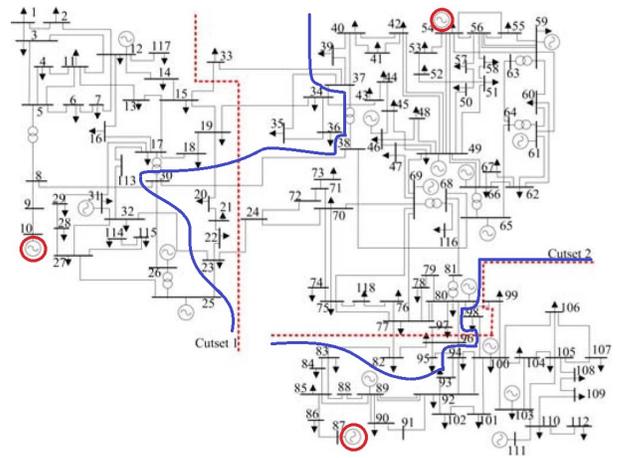

Fig. 2: IEEE 118-bus system topology. The circled generators at buses {10,54,87} are reference generators identified by Algorithm 2. The solid curves represent the cutset selected by Algorithm 1 with parameter $\xi = 10^{-6}$. The dotted lines represent the cutset selected by the spectral clustering method.

TABLE III: Comparison of our proposed islanding ($\xi = 10^{-6}$) and spectral clustering islanding results on the 118-bus system.

| Method | $J(S)$ | Load-generation Imbalance $\sqrt{f(S)}$ | Generator Non-coherency $\overline{H}$ |
|---|---|---|---|
| Proposed Islanding | 0.0524 | 9.7 MW | 1.4515 |
| Spectral Clustering | 0.0527 | 22.1 MW | 1.4515 |

## VII. Conclusion

In this paper, we studied the problem of selecting a set of transmission lines to trip in order to partition an unstable power system into internally stable islands. Our formulation captures the generator coherency and load-generation imbalance jointly in one objective function. We exploited the weak submodularity of the formulation and proposed a greedy islanding algorithm with provable optimality guarantees, which is currently not available in the existing literature. We also studied the reference generator selection and proposed a submodular maximization approach. Our proposed islanding algorithms were validated on the IEEE 39-bus and 118-bus power systems, with a comparison to the state-of-the-art spectral clustering-based islanding method. We observed that the proposed islanding approach can significantly reduce the power imbalance compared to the state of the art, with the same or slightly higher generator non-coherency.